\documentclass[12pt]{article}
\usepackage[mathscr]{eucal}
\usepackage{amssymb,amsmath}
\font\cmc=cmcsc10  scaled \magstep2

\newcommand\N{\mathbb{N}}

\newcommand\vk{\vskip}
\newcommand\hk{\hskip}

\newcommand\al{\alpha}

\newcommand\iy{\infty}

\newcommand\la{\lambda}

\newcommand\sg{\sigma}
\newcommand\el{\ell}

\newcommand\varep{\varepsilon}
\newcommand\rg{\rightarrow}

\newcommand\Om{\Omega}

\newcommand\De{\Delta}

\newcommand\ov{\overset}
\newcommand\und{\underset}
\newcommand\no{\noindent}

\newcommand\ovl{\overline}
\newcommand\for{\forall}
\newcommand\col{\colon\hk-.5em}

\newtheorem{Proof.}{\it Proof.}

\begin{document}
\vbox to .5truecm{}

\begin{center}
\cmc A Note on Self-extremal Sets in {\mathversion{bold} $L_p(\Om)$} Spaces
\end{center}
\vk.3cm

\begin{center}
by Viet Nguyen-Khac
\footnote{Institute of Mathematics, 18 Hoang Quoc Viet, 
10307 Hanoi, Vietnam}  
\ \& Khiem Nguyen-Van
\footnote{Department of Mathematics \& Informatics, Hanoi Univ. of Education, Cau Giay dist., Hanoi, Vietnam.}
\end{center}
\begin{center} 
Hanoi Institute of Mathematics \& Hanoi University of Education
\end{center}
\vk.1cm
\begin{center}
1991 Mathematics Subject Classification. Primary 46B20, 46E30.
\end{center}

\normalsize
\begin{abstract} {We give a geometric characterization of self-extremal sets in $L_p(\Om)$ spaces that partially extends 
our previous results to the case of\ $L_p$\ spaces.}
\end{abstract}

\vk.5cm					
In \cite{Nkk1}, \cite{Nkk2} we introduced the notion of (self-) extremal sets of a Banach space\ $(X,\|\cdot\|)$.\ For a 
non-empty bounded subset\ $A$\ of\ $X$\ we denote by\ $r(A)$\ the relative Chebyshev radius of\ $A$\ with respect to the 
closed convex hull\ $\ovl{co} A$\ of\ $A$,\ {\it i.e.}\ $r(A)\col=\und{y\in \ovl{co} A}{\text {inf}}\ \und{x\in A}{\text {sup}}\ 
\|x-y\|$.\ The self-Jung constant of\ $X$\ is defined by\ $J_s(X)\col=$ sup $\{r(A):\ A\subset X,\ \ {\text {with\ \ diam}}(A)=
1\}$.\ If in this definition we replace\ $r(A)$\ by the relative Chebyshev radius\ $r_X(A)$\ of\ $A$\ with respect to the 
whole\ $X$,\ we get the Jung constant\ $J(X)$\ of\ $X$.\ Recall that a bounded subset\ $A$\ of\ $X$\ consisting of at least 
two points is said to be extremal (resp. self-extremal) if\ $r_X(A)=J(X)\ d(A)$\ (resp.\ $r(A)=J_s(X)\ d(A)$).
\vk.3cm
Let\ $(\Om,\mu)$\ be a\ $\sg$-finite measure space such that\ $L_p(\Om)$\ is infinite-dimentional. The Jung and 
self-Jung constants of\ $L_p(\Om)\ \ (1\le p<\iy)$\ were determined in \cite{DB}, \cite{IvP}, \cite{Pich1}, \cite{Pich2}: 
   
  $$J(L_p(\Om))=J_s(L_p(\Om))=\max\ \{ 2^{\frac{1}{p}-1},\ 2^{-\frac{1}{p}}\}.$$

\vk.3cm 
\textbf{Theorem\ 1.1.}\quad {\it Let\ $(\Om,\mu)$\ be a\ $\sg$-finite measure space such that\ $L_p(\Om)$\ is 
infinite-dimentional\ $(1<p<\iy)$.\ If\ $A$\ is a self-extremal set in\ $L_p(\Om)$,\ then\ $\al(A)=d(A)$.}
\vk.3cm
Here\ $\al(A)\col=\inf\ \{\varep>0:\ A$\ can be covered by finitely many sets of diameter\ $\le\varep\}$\ - the 
Kuratowski measure of non-compactness of\ $A$.
\vk.3cm
Before proving our theorem we need the following results which for convenience we reformulate in the form of 
Lemmas 1.2 and 1.3 below.
\vk.3cm
\textbf{Lemma\ 1.2} (\cite{DB}, Theorem 1.1).\quad {\it Let\ $X$\ be a reflexive strictly convex Banach space and\ $A$\ 
a finite subset of\ $X$.\ Then there exists a subset\ $B\subset A$\ such that
\vk.2cm
{\rm (i)}\ \ $r(B)\ge r(A)$;
\vk.2cm
{\rm (ii)}\ \ $\|x-b\|=r(B)$\ for every $x\in B$,\ where\ $b$\ is the relative Chebyshev center of\ $B$,\ {\it i.e.}\ $b\in 
\ovl{co} B$\ and\ $\und{x\in B}{\text {sup}}\ \|x-b\|=r(B)$.}
\vk.3cm
\textbf{Lemma\ 1.3} (\cite{WW}, Theorem 15.1).\quad {\it Let\ $(\Om,\mu)$\ be a\ $\sg$-finite measure space,\ $1<p<
\iy,\ x_1, \ldots, x_n$\ vectors in\ $L_p(\Om)$\ and\ $t_1, \ldots, t_n$\ non-negative numbers such that\ $\und{i=1}
{\ov{n}{\sum}}\ t_i=1.$\ The following inequality holds
  
  $$2 \sum_{i=1}^n\ t_i \Big\|x_i-\sum_{j=1}^n\ t_j x_j\Big\|^\al\le \sum_{i, j=1}^n\ t_i t_j \| x_i-x_j\|^\al,$$

\no where
  
  $$\al=\begin{cases}  q\col=\dfrac{p}{p-1}\ \ & {\text if}\ 1<p\le 2\\
                                    p\ \ & {\text if}\ p\ge 2\end{cases}\eqno{(1)}$$   }
\vk.3cm
{\it Proof of Theorem 1.1}.\quad We may assume that\ $A$\ is closed convex and\ $r(A)=1$.\ For each integer\ 
$n\ge 2$\ we have\ 
  
  $$\und{x\in A}{\bigcap}\ B(x,1-\dfrac{1}{n})\ \cap\ A=\emptyset,$$

\no where\ $B(x,r)$\ denotes the closed ball centered at\ $x$\ with radius\ $r$\ which is weakly compact since\ $L_p(\Om)$\ 
is reflexive. Hence there exist\  $x_{t_{n-1}+1},\ x_{t_{n-1}+2},\ \ldots,\ x_{t_n}$\ in\ $A$\ (with convention\ $t_1=0$) 
such that 
  
  $$\und{i=t_{n-1}+1}{\ov{t_n}{\bigcap}}\ B(x_i,1-\dfrac{1}{n})\ \cap\ A=\emptyset$$

Set\ $A_n\col=\{x_{t_{n-1}+1},\ x_{t_{n-1}+2},\ \ldots,\ x_{t_n}\}$.\ By Lemma 1.2 there exists a subset\ $B_n=
\{y_{s_{n-1}+1},\ y_{s_{n-1}+2},\ \ldots,\ y_{s_n}\}$ of\ $A_n$\ satisfying properties (i)-(ii) of the lemma.\ Let us denote 
the relative Chebyshev center of\ $B_n$\ by\ $c_n$\ and let\ $r_n\col=r(B_n)$.\ By what we said above we have\ $r_n>
1-\dfrac{1}{n}$\ and\ $\|y_i-c_n\|=r_n$\ for every\ $i\in I_n\col=\{s_{n-1}+1,\ s_{n-1}+2,\ \ldots,\ s_n\}$. Since\ $B_n$\ 
is a finite set there exist non-negative numbers\ $\al_{s_{n-1}+1},\ \al_{s_{n-1}+2},\ \ldots,\ \al_{s_n}$\ with\ 
$\und{i\in I_n}{\sum}\ \al_i=1$\ such that\ $c_n=\und{i\in I_n}{\sum}\ \al_i y_i$.\ Applying Lemma 1.3 one gets    
  
  $$2 r_n^\al=2 \und{i\in I_n}{\sum}\ \al_i \Big\|y_i-\und{j\in I_n}{\sum}\ \al_j y_j\Big\|^\al\le \und{i, j\in I_n}{\sum}\ 
  \al_i \al_j\| y_i-y_j\|^\al\eqno{(2)}$$

\no where\ $\al$\ as in (1). 
\vk.2cm
Setting\ $B_\iy\col=\{y_{s_{n-1}+1},\ y_{s_{n-1}+2},\ \ldots,\ y_{s_n}\}_{n=2}^\iy$\ we claim that\ $\al(B_\iy)=d(A)$.\ 
Suppose on the contrary\ $\al(A_\iy)<d(A)$,\ so there exists\ $\varep_0\in (0,d(A))$\ satisfying\ $\al(B_\iy)\le 
d(A)-\varep_0$,\ and subsets\ $D_1,\ D_2,\ \cdots,\ D_m$\ of\ $L_p(\Om)$\ with\ $d(D_i)\le d(A)-\varep_0$\ for 
every\ $i=1,\ 2,\ \ldots,\ m$,\ such that\ $B_\iy\subset\und{i=1}{\ov{m}{\bigcup}}\ D_i$.\ Then one can find at least one 
set among\ $D_1,\ D_2,\ \cdots,\ D_m$,\ say\ $D_1$\ with the property that there are infinitely many\ $n$\ satisfying
  
  $$\sum_{i\in J_n}\ \al_i\ge \dfrac{1}{m}\eqno{(3)}$$

\no where
  
  $$J_n\col=\{i\in I_n\colon\ \ y_i\in D_1\}.$$

\vk.2cm
In view of (2) and taking into account\ $\big( d(A)\big)^\al=\Big(\dfrac{1}{J_s\big( L_p(\Om)\big)}\Big)^\al=2$\ we have 
for all\ $n$\ satisfying\ (3) 

  $$\begin{aligned} 2.r_n^\al & \le \und{i, j\in I_n}{\sum}\ \al_i.\al_j.\|y_i-y_j\|^\al\\
  & \le \big( d(A)-\varep_0 \big)^\al.\Big(\und{i, j\in J_n}{\sum}\ \al_i.\al_j\Big)+\big( d(A)\big)^\al.\Big( 1-
  \und{i, j\in J_n}{\sum}\ \al_i.\al_j\Big)\\
  & \le 2-\Big[ \big( d(A)\big)^\al-\big( d(A)-\varep_0 \big)^\al\Big] .\dfrac{1}{m^2}
  \end{aligned}\eqno{(4)}$$
\vk.2cm
On the other hand, obviously\ $1-\dfrac{1}{n}< r_n\le 1$,\ therefore\ $\und{n\to\iy}{\text {lim}}\ r_n=1$.\ We get a 
contradiction with (4), since there are infinitely many\ $n$\ satisfying (3).  
\vk.2cm
One concludes that\ $\al(B_\iy)=d(A)$,\ and hence\ $\al(A)=d(A)$.
\vk.2cm
The proof of Theorem 1.1 is complete.
\vk.5cm
As an immediate consequence one obtains an extension of Gulevich's result for\ $L_p(\Om)$\ spaces. 
\vk.5cm
\textbf{Corollary\ 1.4}\ ({\it cf.} \cite{Gul}).\quad {\it Let\ $(\Om,\mu)$\ be a\ $\sg$-finite measure space such that\ 
$L_p(\Om)$\ is infinite-dimentional\ $(1<p<\iy)$.\ Assume\ $A$  a relatively compact set in\ $\L_p(\Om)$\ with\ 
$d(A)>0$.\ Then

  $$\begin{aligned} & r(A)<\dfrac{1}{\sqrt[q]{2}}\cdot d(A),\qquad 1<p\le 2,\ q\col=\dfrac{p}{p-1};\\
    & r(A)<\dfrac{1}{\sqrt[p]{2}}\cdot d(A),\qquad 2<p<\iy.
  \end{aligned}$$}
\vk.3cm
\textbf{Theorem\ 1.5.}\quad {\it Under the assumptions of Theorem 1.1, for every\ $\varep\in (0,d(A))$,\ every positive 
integer\ $m$,\ there exists an\ $m$-simplex\ $\De(\varep,m)$\ with vertices in\ $A$\ such that each edge of\ 
$\De(\varep,m)$\ has length not less than\ $d(A)-\varep$.}
\vk.2cm
{\it Proof}.\quad We shall assume\ $A$\ is closed convex, and\ $r(A)=1$.\ From the proof of Theorem 1.1 we derived a 
sequence\ $\{y_{s_{n-1}+1},\ y_{s_{n-1}+2},\ \ldots,\ y_{s_n}\}^\iy_{n=2}$\ in\ $A$,\ and a sequence of positive 
numbers\ $\{\al_{s_{n-1}+1},\ \al_{s_{n1}+2},\ \ldots,\ \al_{s_n}\}^\iy_{n=2}$\ (with convention\ $s_1=0$)\ such that
  $$2.r^\al_n\le \und{i, j\in I_n}{\sum}\ \al_i.\al_j.\|y_i-y_j\|^\al,\ \ \und{i\in I_n}{\sum}\ \al_i=1,$$

\no where\ $r_n\in (1-\frac{1}{n},1],\ \al$\ as in (1), and\ $I_n\col=\{s_{n-1}+1,\ s_{n-1}+2,\ \ldots,\ s_n\}$.   
\vk.2cm
We denote by

  \begin{align} & T_{nj}\col=\und{i\in I_n}{\sum}\ \al_i.\|y_i-y_j\|^\al,\notag\\
  & S_n\col=\big\{ j\in I_n:\ \ T_{nj}\ge 2.r^\al_n .\big( 1-\sqrt{1-r^\al_n} \big) \big\},\notag\\
  & S_n(y_j)\col=\bigg\{ i\in I_n:\ \ \|y_i-y_j\|^\al\ge 2.\Big( 1-\dfrac{1}{\sqrt[4]{n}}\Big)\bigg\},\ j\in S_n,\notag\\
  & \hat{S}_n(y_j)\col=\big\{ y_i:\ \ i\in S_n(y_j) \big\},\ j\in S_n,\notag\\
  & \la_n\col=\und{i\in I_n\setminus S_n}{\sum}\ \al_i=1-\und{i\in S_n}{\sum}\ \al_i. \notag
  \end{align}
\vk.2cm
Based on the method in proving Theorem 3.4 of \cite{Nkk2} one can proceed furthermore as follows. We have
  
  $$\begin{aligned} 2\ r^\al_n & \le \und{i, j\in I_n}{\sum}\ \al_i\ \al_j\ \|y_i-y_j\|^\al\\
  & =\und{j\in S_n}{\sum}\ \al_j\ \und{i\in I_n}{\sum}\ \al_i\ \|y_i-y_j\|^\al+\und{j\in I_n\setminus S_n}{\sum}\ \al_j\ 
   \und{i\in I_n}{\sum}\ \al_i\ \|y_i-y_j\|^\al\\
  & \le 2 \und{j\in S_n}{\sum}\ \al_j+2\ r^\al_n\ \big( 1-\sqrt{1-r^\al_n} \big)\und{j\in I_n\setminus S_n}{\sum}\ \al_j\\
  & =2-2\ \la_n\ \big( 1-r^\al_n+r^\al_n\ \sqrt{1-r^\al_n} \big)\\
  & \le 2-2\ \la_n\ \sqrt{1-r^\al_n}.
 \end{aligned}$$

Hence $\la_n\le \sqrt{1-r^\al_n}\rg 0$,\ as\ $n\rg\iy$. Thus $\und{n\to\iy}{\text {lim}}\ \Big(\und{i\in S_n}{\sum}\ 
\al_i\Big)=\und{n\to\iy}{\text {lim}}\ (1-\la_n)=1.$
\vk.2cm
On the other hand

 $$2\ r^\al_n\le \und{i, j\in I_n}{\sum}\ \al_i\ \al_j\ \|y_i-y_j\|^\al \le 2\ \bigg( 1-\Big(\und{i\in I_n}{\sum}\ \al_i^2\Big)\bigg)
 \le 2\ (1-\al_i^2)$$

\no for every\ $i\in I_n$.\ Therefore\ $\al_i\le \sqrt{1-r^\al_n}\rg 0$\ as\ $n\rg\iy$.\ One concludes that the cardinality\ 
$\big| S_n\big|$\ of\ $S_n$\ tends to\ $\iy$\ as\ $n\rg\iy.$\ In a similar manner as in the proof of Theorem 3.4 of 
\cite{Nkk2} for every\ $\varep\in (0,d(A))$\ and a given positive integer\ $m$\ we choose\ $n$\ sufficiently large satisfying
 $$\big| S_n\big|>m,\quad\ \dfrac{2\ \al\ m}{\sqrt[4]{n}}<1,\quad\ 2 \Big( 1-\dfrac{1}{\sqrt[4]{n}}\Big)\ge \big( d(A)-
\varep\big)^\al$$ 
such that for every\ $1\le k\le m$\ and every choice of\ $i_1,\ i_2,\ \ldots,\ i_k\in S_n$\ we have
  
   $$\und{\nu=1}{\ov{k}{\bigcap}}\ \hat{S}_n(y_{i_\nu})\ne\emptyset.$$
\vk.2cm
With\ $m$\ and\ $n$\ as above and a fixed\ $j\in S_n$,\ setting\ $z_1\col=y_j$,\ we take consecutively\ 
$z_2\in \hat{S}_n(z_1),\ z_3\in \hat{S}_n(z_1)\cap \hat{S}_n(z_2),\ \ldots,\ z_{m+1}\in \und{k=1}{\ov{m}{\bigcap}}\ 
\hat{S}_n(z_k)$. One sees that

  $$\|z_i-z_j\|^\al\ge 2 \Big( 1-\dfrac{1}{\sqrt[4]{n}}\Big)\ge \big( d(A)-\varep\big)^\al.$$

\no for all\ $i\ne j$\ in\ $\{1, 2, \ldots, m+1\}$, with\ $n$\ sufficiently large. We obtain an\ $m$-simplex\ formed by\ 
$z_1,\ z_2,\ \ldots,\ z_{m+1}$,\ whose edges have length not less than\ $d(A)-\varep$,\ as claimed.
\vk.2cm
The proof of Theorem 1.5 is complete.
\vk.3cm
\textbf{Remark\ 1.6.}\quad (i)\ Since for\ $L_p(\Om)$\ spaces\ $J_s=J$\ the extremal sets in\ $L_p(\Om)$\ are also 
self-extremal. Thus we obtain a characterization for extremal sets in\ $L_p(\Om)$ via Theorem 1.5 above.
\vk.2cm
(ii)\ In particular\ $\Om=\N,\ \mu(A)\col=$ card $(A),\ A\subset\N$ leads to the\ $\ell_p$\ space case (Theorem 
3.4 of \cite{Nkk2}).
\vk.3cm
\textbf{Example\ 1.7.}\quad (i)\ Let\ $p\ge 2$, consider a sequence\ $\{\Om_n\}_{i=1}^\iy$\ consisting of measurable 
subsets of\ $\Om$\ such that
  
$$0<\mu(\Om_i)<\iy,\ i=1, 2, \ldots;\quad \Om_i\cap\Om_j=\emptyset,\ \for\ i\ne j;\quad \bigcup_{i=1}^\iy\ \Om_i=
\Om.$$
\vk.2cm
Let\ $\chi_{\Om_i}$\ denote the characteristic function of\ $\Om_i$,\ and set

  $$A\col=\{ f_i\}_{i=1}^\iy,\ \ {\text {where}}\ \ f_i\col=\dfrac{\chi_{\Om_i}}{\mu(\Om_i)}.$$

\no One can check easily that\ $r(A)=1,\ d(A)=2^{\frac{1}{p}}$,\ hence\ $A$\ is a self-extremal set in\ $L_p(\Om)$.
\vk.3cm
(ii)\ In the case\ $1<p<2$\ we set\ $B\col=\{ r_i\}_{i=0}^\iy,$\ where\ $\{r_i\}_{i=0}^\iy$\ is the sequence of 
Rademacher functions in\ $L_p[0,1]$.\ If\ $r\in co \{r_0, r_1, \ldots, r_n\}$\ and\ $k\ge n+1$\ then it is easy to see that\ 
$d(B)=2^{1-\frac{1}{p}}$\ and 

  $$\|r-r_k\|_p\col=\Big(\int_0^1\ |r-r_k|^p\ d\mu\Big)^{\frac{1}{p}}\ge \Big |\int_0^1\ (r-r_k) r_k\ d\mu
\Big |=1,$$

\no hence\ $r(B)=1$.\ Thus\ $B$\ is a self-extremal set in\ $L_p[0,1]$\ with\ $1<p<2$.  This is in contrast to the\ 
$\ell_p$\ case (\cite{Nkk2}), where we conjectured that there are no (self)-extremal sets in\ $\ell_p$,\ spaces with\ 
$1<p<2$.
 
\vk.5cm

\end{document}